\documentclass[a4paper,14pt]{article}
\usepackage[cp1251]{inputenc}
\usepackage{float}
\usepackage{cite}
\usepackage[english]{babel}
\usepackage{amssymb,amsmath,amsfonts}
\usepackage{enumerate}
\usepackage[dvips]{graphicx}
\usepackage{setspace}
\usepackage{fancyhdr}

\usepackage{xcolor}
\usepackage{ulem}  

\begin{document}

\begin{center}
{\textbf{Finite groups with Hall 2-maximal subgroups}}\\
\vspace{\baselineskip}

M.\,N.~Konovalova$^1$, V.\,S.~Monakhov$^2$\\

\textit{$^1$Russian Presidential Academy of National Economy
and Public Administration, Bryansk,} e-mail: msafe83@mail.ru;\\

\textit{$^2$F. Scorina Gomel State University,}\\e-mail: victor.monakhov@gmail.com
\vspace{\baselineskip}

\end{center}

\vspace{\baselineskip}

All groups in this paper are finite.
The notation $M\le G$ ($M<G$, $M\lessdot \,G$) means that~$M$
is a subgroup (proper subgroup, maximal subgroup, respectively)
of a group~$G$.
A subgroup $H$ of a group $G$ is called a Hall subgroup
if $|H|$ and $|G:H|$ are coprime.

V.\,S.~Monakhov~\cite{Pi1}, N.\,V.~Maslova and
D.\,O.~Revin~\cite{Pi2} studied groups with Hall maximal subgroups.

A subgroup $H$ of a group $G$ is said to be a $2$-maximal subgroup
if there is a maximal subgroup $M$ of $G$ such that $H$ is a
maximal subgroup of $M$.

The group structure substantially depends on
the properties of its 2-maximal subgroups.
Huppert, Suzuki, Yanko, V.\,A.~Belonogov got the earliest results.
Huppert~\cite{Pi3} established the supersolubility
of groups in which all 2-maximal subgroups are normal.
Suzuki~\cite{Pi4} and Yanko~\cite{Pi5} proved that
an unsolvable group with nilpotent 2-maximal subgroups
is isomorphic to~$PSL(2,5)$ or~$SL(2,5)$.
V.\,A.~Belonogov~\cite{Pi6} described finite solvable groups
with nilpotent 2-maximal subgroups.
In the current century, this direction
was investigated in many papers,
see references in~\cite{Pi7}--\cite{Pi10}.

In this paper,
we study a group in which every 2-maximal subgroup
is a Hall subgroup.
Based on the results of
N.\,V.~Maslova and D.\,O.~Revin~\cite{Pi2},
we obtain the solvability of a group,
and then we apply the description from~\cite{Pi1}.

We need the following definitions and notation.

\medskip

{\sl A Gasch\"{u}tz subgroup} of a group $G$ is a subgroup $W$
such that

$(1)$ $W$ is supersolvable;

$(2)$ if $W\leq A<B\leq G$, then $|B:A|$~ is not a prime.

\medskip

In a solvable group, Gasch\"{u}tz subgroups exist and are
conjugate~\cite[5.29]{Pi11}.

\medskip

{\sl The supersolvable residual} $G^\mathfrak{U}$
of a group $G$ is called the smallest normal subgroup in~$G$
whose quotient is supersolvable.

\medskip

Let $G$~ be a group. Then

$|G|$~is the order of~$G$;

$\pi (G)$~is the set of all prime divisors of $|G|$;

$G_p$ is a Sylow $p$-subgroup of~$G$;

$\sigma (G)=\{r\in \pi (G)\, : \, |G_r|=r\}$;

$\tau (G)=\{r\in \pi (G)\, : \, |G_r|>r\}$.

It is clear that~$\pi (G)=\sigma (G)\cup \tau (G)$
and~$\sigma (G)\cap \tau (G)$~is an empty set.

\medskip

{\bf Theorem.} {\sl Let $G$~be a finite non-primary group in which each
2-maximal subgroup is a Hall subgroup.
If $G$~is a supersolvable group, then
$|G|$ is square-free, that is
$\pi (G)=\sigma (G)$. Let $G$~be a non-supersolvable group.
Then the following statements hold.

$(1)$ $G$ has a Sylow tower and
each maximal subgroup of $G$ is a Hall subgroup.

$(2)$~$|\sigma (G)|\ge 2$ and $G_{\sigma (G)}\le W$, where $W$ is
a Gasch\"{u}tz subgroup of $G$.

$(3)$~$|\tau (G)|\ge 1$ and~$G_{\tau (G)}$ is
the supersolvable residual of $G$.}

\medskip

{\bf Proof.}
Denote $\sigma (G)=\sigma$ and~$\tau(G)=\tau$.

Let~$M$~be a maximal subgroup in~$G$ and~$H$~be a maximal subgroup
in~$M$. Then~$H$~is a 2-maximal subgroup of~$G$ and, by
the hypothesis,~$H$ is a Hall subgroup of~$G$.
Assume that~$M$ is not a Hall subgroup.
Then there is~$p\in \pi (G)$ such that
$p$ divides~$|M|$ and~$p$ divides~$|G:M|$.
If $M$ is a $p$-subgroup, then~$H=1$ and~$|G|=p^2$, that is~$|\pi (G)|=1$,
a contradiction to the hypothesis.
Consequently,~$M$~is not a $p$-subgroup,
and we can choose~$H$
so that $H$ is a maximal subgroup of $M$
and $H$ contains a Sylow $p$-subgroup of~$M$.
Then $p$ divides $|H|$ and $p$
divides~$|G:H|=|G:M||M:H|$, a contradiction.
Therefore, the assumption is false,
and each maximal subgroup in $G$ is a Hall subgroup.

If $G$~is supersolvable group, then
$|G|$ is square-free~\cite[Corollary 3]{Pi1}.

Assume that $G$ is an unsolvable group.
By the result of N.\,V.~Maslova and
D.\,O.~Revin~\cite[Corollary 1]{Pi2}, the
solvable radical $S(G)$ has a Sylow tower,
and $G/S(G)$ is either trivial
or isomorphic to one of the groups $PSL(2,7)$, $PSL(2,11)$
or $PSL(5,2)$. Maximal subgroups
of these groups are known.

If $G/S(G)$ is isomorphic to $PSL(2,7)$,
then we choose $M$ and~$H$ in $G$ such that
$$
S(G)\le H\lessdot \, M\lessdot \,G, \
M/S(G)\cong S_4, \ H/S(G)\cong A_4.
$$
Hence $H$ is a 2-maximal subgroup and is not a Hall subgroup in~$G$,
since
$$
|H|=2^2\cdot 3\cdot |S(G)|, \ |G:H|=7\cdot 2\cdot |H|.
$$
If $G/S(G)$ is isomorphic to $PSL(2,11)$,
then we choose $M$ and~$H$ in $G$ such that
$$
S(G)\le H\lessdot \, M\lessdot \,G, \
M/S(G)\cong D_{12}, \ H/S(G)\cong D_6.
$$
So $H$ is a 2-maximal subgroup and is not a Hall subgroup in~$G$,
since
$$
|H|=2\cdot 3\cdot |S(G)|, \ |G:H|=5\cdot 11\cdot 2\cdot |H|.
$$
If $G/S(G)$ is isomorphic to $PSL(5,2)$,
then we choose $M$ and~$H$ in $G$ such that
$$
S(G)\le H\lessdot \, M\lessdot \,G, \
M/S(G)\cong 2^6:(S_3\times PSL(3,2)), \
H/S(G)\cong 2^6:(Z_3\times PSL(3,2)).
$$
Hence $H$ is a 2-maximal subgroup and is not a Hall subgroup in~$G$,
since
$$
|H|=2^{9}\cdot 3^2\cdot 7\cdot |S(G)|, \
|G:H|=2\cdot 5\cdot 31\cdot |S(G)|.
$$
Therefore we conclude that $G/S(G)$ is trivial,
and $G$ is solvable.

Thus $G$ is solvable
and all its maximal subgroups are Hall subgroups.
By~\cite[Corollary 2]{Pi1}, $G$ has a Sylow tower
and each Sylow subgroup in~$G$ is an elementary Abelian subgroup.

Assume that $K\lessdot \, M\lessdot \,G$,
$K$ is normal in~$M$, $M$ is normal in~$G$. Then~$|G:M|=p$, $|M:K|=q$, $p$ and~$q$
are  primes.
Since~$M$ and~$K$~are Hall subgroups,  we obtain $p\ne q$,
$p=|G_p|\in \sigma$, $q=|G_q|\in \sigma$, and~$|\sigma|\ge 2$.
Let~$W$  be a Gasch\"{u}tz subgroup
of $G$.
By definition, in the subgroups chain
$$
W=W_0\lessdot \, W_1\lessdot \, \ldots W_{m-1}\lessdot \, W_m=G
$$
all indices
$|W_{i+1}:W_i|$~are not primes.
Since all these indices are primary,
we have $\pi (G:W)\subseteq \tau$
and~$G_{\sigma}\le W$ for some Hall
$\sigma$-subgroup~$G_{\sigma}$ of~$G$.

Let~$H=G^\frak U$ be the
supersolvable residual of $G$.
Since $G$ is not supersolvable,
we get $H\ne 1$ and~$|\tau|\ge 1$.
In $G/H$, all 2-maximal subgroups
are Hall subgroups,
therefore~$\pi (G:H)\subseteq \sigma$ and~$G_\tau \le H$.
By~\cite[Corollary 1]{Pi1}, $H$ is a Hall subgroup.
If~$F(G)$ is a $\sigma$-subgroup,
then~$F(G)$ is cyclic and $G$ is supersolvable
in view of \cite[Corollary 3]{Pi1},
a contradiction. Therefore $F(G)$ is not a
$\sigma$-subgroup, and there is
a nontrivial normal in~$G$ $r$-subgroup~$R$ for~$r\in \tau$.
It is clear that~$R\le G_\tau \le H$.
By~\cite[Corollary 2]{Pi1}, $R$ is a Sylow subgroup of~$G$
and minimal normal subgroup in~$G$.
If $G/R$ is supersolvable, then~$H=R$ and~$H=G_\tau$.
If~$G/R$ is not supersolvable, then by the induction hypothesis,
$$
(G/R)^\frak U=G^\frak U/R=(G/R)_\tau =G_\tau /R, \ G^\frak U=H=G_\tau.
$$
The theorem is proved.

\medskip

{\bf Example.}
The general linear group $GL_{2}(29)$ has a cyclic
subgroup $Z_ {15}$ of order $15$, which acts irreducibly on
an elementary abelian group $E_{29^2}$ of order $29^2$.
In the semidirect product $E_{29^2}\rtimes Z_{15}$,
all maximal subgroups and all 2-maximal subgroups
are Hall subgroups.

\end{document}